\documentclass[11pt,letterpaper,twoside,spanish]{amsart}

\usepackage{multicol}
\usepackage{amssymb}
\usepackage{amsmath}
\usepackage{amsbsy}
\usepackage{amsthm}
\usepackage[latin1]{inputenc}

\usepackage[spanish]{babel}

\newtheorem{thm}{Theorem}
\newtheorem{cor}[]{Corollary}
\newtheorem{lem}[]{Lemma}

\theoremstyle{remark}
\newtheorem{rem}[]{Remark}

\newcommand{\sop}{\text{supp}}

\newcommand{\R}{\mathbb{R}}

\newcommand{\C}{\mathbb{C}}

\begin{document}


\title{Fourier-Pad{\'e} approximants for Angelesco systems}

\author{M. Bello-Hern{\'a}ndez}
\address{Universidad de La Rioja \\ Dpto. Matem{\'a}ticas y Computaci{\'o}n}

\author{G. L{\'o}pez-Lagimasino}
\address{Universidad Carlos III de Madrid\\ Dpto. de Matem{\'a}ticas
Aplicada}

\author{J. M{\'\i}nguez-Ceniceros}
\address{Universidad de La Rioja\\ Dpto. Matem{\'a}ticas y Computaci{\'o}n}

\maketitle


\vspace{1cm}


\section{Introduction}

In this paper we study linear and non-linear Fourier-Pad\'e
approximation for Angelesco systems of functions. This
construction is similar to that of Hermite-Pad\'e approximation.
Instead of considering power series expansions of the functions in
the system, we take their expansion in a series of orthogonal
polynomials.

In \cite{sue} and \cite{sue2}, S. P. Suetin obtained convergence
results for rows of Fourier-Pad\'e approximation extending to this
setting classical results of the theory of Pad\'e approximation.

Diagonal sequences of Fourier-Pad\'e approximation were studied by
A. A. Gonchar, E. A. Rakhmanov, and S. P. Suetin in
\cite{Gon-Rak-Sue} when the function to be approximated is of
Markov type; that is, the Cauchy transform of a measure supported
on the real line. They give   the rate of convergence of diagonal
sequences of linear and non-linear Pad\'e approximants in terms of
the equilibrium measures of a related   potential theoretic
problem. We generalize those  results   to the case when  a system
of Markov functions is given defined by measures whose supports do
not intersect.

Let $\mathcal{M}(\Delta)$ denote the class of all finite, Borel
measures with compact support consisting of an infinite set of
points contained in an interval $\Delta$ of the real line
$\mathbb{R}$. Given $\sigma \in \mathcal{M}(\Delta)$, let
\[
\widehat{\sigma}(z)=\int\frac{d\sigma(x)}{z-x}
\]
be the associated  \emph{Markov} function. Let $\Delta_k,
k=1,\ldots,m,$ be  intervals of the real line such that
\[ \Delta_k \cap \Delta_j = \emptyset\,, \qquad  k \neq j \,,
\]
and $\sigma_k \in \mathcal{M}(\Delta_k), k=1,\ldots,m.$ We say
that   $\sigma = (\sigma_1,\ldots,\sigma_m)$ forms an Angelesco
system of measures and
$(\widehat{\sigma}_1,\ldots,\widehat{\sigma}_m)$ is the associated
Angelesco system of functions.

Let $\sigma_0 \in \mathcal{M}(\Delta_0).$ Likewise, we will assume
that
\[ \Delta_0 \cap \Delta_k = \emptyset\,, \qquad k=1,\ldots,m\,.
\]
Consider the sequence $\{\ell_n\}, n \in \mathbb{Z}_+ =
\{0,1,2,\ldots\},$ of orthonormal polynomials with respect to
$\sigma_0$ with positive leading coefficient. Take a multi-index
${\bf n} = (n_1,\ldots,n_m) \in \mathbb{Z}_+^m$. Set
\[ |{\bf n}| = n_1 + \cdots + n_m \,.
\]

Let $Q_{\bf n},P_{{\bf n},1},\ldots,P_{{\bf n},m},$ be polynomials
such that:
\begin{itemize}
\item[i)] $ \deg Q_{\bf n} \leq |{\bf n}|, Q_{\bf n} \not\equiv 0\,,
\deg P_{{\bf n},j} \leq |{\bf n}|-1, j=1,\ldots,m\,.$
\item[ii)] For each $j = 1,\ldots,m\,,$ and $k = 0,\ldots, |{\bf
n}| +n_j -1$
\[
c_k(Q_{\bf n} \widehat{\sigma}_j - P_{{\bf n},j}) = \int (Q_{\bf
n} \widehat{\sigma}_j - P_{{\bf n},j})(x) \ell_k(x) d \sigma_0(x)
= 0\,.
\]
\end{itemize}
Notice that
\begin{equation} \label{losP}
 P_{{\bf n},j}(z) = \sum_{i=0}^{|{\bf
n}|-1}c_i(Q_{\bf n} \widehat{\sigma}_j)\ell_i(z)\,.
\end{equation}
The $|{\bf n}| +1$ coefficients of $Q_{\bf n}$ satisfy a
homogeneous linear system of $|{\bf n}|$ equations given by
\[ c_k(Q_{\bf n}\widehat{\sigma}_j) = 0\,, \qquad j = 1,\ldots,m\,,\qquad k
= |{\bf n}|,\ldots,|{\bf n}| + n_j -1\,.
\]
Therefore, a non-trivial solution is guaranteed.

In Section \ref{teo-pot} we will prove that every solution to
i)-ii) has $\deg Q_{\bf n} = |{\bf n}|$. This being the case,
$(Q_{\bf n},P_{{\bf n},1},\ldots,P_{{\bf n},m})$ is uniquely
determined up to a constant factor. In fact, let us assume that
$(Q_{\bf n},P_{{\bf n},1},\ldots,P_{{\bf n},m}),$ and
$(\widetilde{Q}_{\bf n},\widetilde{P}_{{\bf
n},1},\ldots,\widetilde{P}_{{\bf n},m}),$ are solutions of i)-ii).
Without loss of generality, we can assume that $Q_{\bf n}$ and
$\widetilde{Q}_{\bf n}$ are monic (with leading coefficient equal
to one). Obviously, if $Q_{\bf n} - \widetilde{Q}_{\bf n}
\not\equiv 0$ then $(Q_{\bf n }- \widetilde{Q}_{\bf n},{P}_{{\bf
n},1}-\widetilde{P}_{{\bf n},1},\ldots,{P}_{{\bf n},m} -
\widetilde{P}_{{\bf n},m})$ is also a solution with $\deg Q_{\bf n
}- \widetilde{Q}_{\bf n} < |{\bf n}|$ which contradicts our
assumption. Hence $Q_{\bf n } \equiv \widetilde{Q}_{\bf n}$ and by
(\ref{losP}) it follows that ${P}_{{\bf n},j} \equiv
\widetilde{P}_{{\bf n},j} , j = 1\ldots,m\,.$

The rational vector function $\left(\frac{P_{{\bf n},1}}{Q_{\bf
n}},\ldots,\frac{P_{{\bf n},m}}{Q_{\bf n}}\right)$ constructed
from any solution of i)-ii) is called the {\bf n}-th linear
Fourier-Pad\'e approximant for the Angelesco system
$(\widehat{\sigma}_1,\ldots,\widehat{\sigma}_m)$ with respect to
$\sigma_0$. We shall see that for all ${\bf n} \in \mathbb{Z}_+^m$
the linear Fourier-Pad\'e approximant of an Angelesco system is
unique.

Non-linear Fourier-Pad\'e approximants are determined as follows.
Given ${\bf n} \in \mathbb{Z}_+^m,$ we must find polynomials
$T_{{\bf n}},\, S_{{\bf n},1},\ldots,S_{{\bf n},m}$ such that
\begin{itemize}
\item[i')]
$\deg T_{\bf n} \le |{\bf n}|, T_{\bf n} \not\equiv 0\,,
\deg(S_{{\bf n},j})\le |{\bf n}|-1,\, j=1,\ldots, m\,.$
\item[ii')] For each $j = 1,\ldots,m\,,$ and $k = 0,\ldots, |{\bf
n}| +n_j -1$
\[
c_k\left( \widehat{\sigma}_j - \frac{S_{{\bf n},j}}{T_{\bf n}}
\right) = \int \left( \widehat{\sigma}_j - \frac{S_{{\bf
n},j}}{T_{\bf n}}\right)(x) \ell_k(x) d \sigma_0(x) = 0\,.
\]
\end{itemize}
This system of equations is non-linear in the coefficients of the
polynomials. We shall prove that for each  ${\bf n} \in
\mathbb{Z}_+^m,$ the system has a solution but we have not been
able to show that it is unique. For any solution of i')-ii'), the
vector rational function $\left(\frac{S_{{\bf n},1}}{T_{\bf
n}},\ldots,\frac{S_{{\bf n},m}}{T_{\bf n}}\right)$ is called an
{\bf n}-th non-linear Fourier-Pad\'e approximant for the Angelesco
system $(\widehat{\sigma}_1,\ldots,\widehat{\sigma}_m)$ with
respect to $\sigma_0$.

In this paper, we obtain the rate of convergence (divergence) of
linear and non-linear Fourier-Pad\'e approximants for Angelesco
systems such that the measures $\sigma_0,\ldots, \sigma_m$ are in
the class $\mbox{\bf  {Reg}}$ of regular measures. For different
equivalent forms of defining regular measures see sections 3.1 to
3.3 in \cite{stto}. In particular, $\sigma_0 \in \mbox{\bf {Reg}}$
if and only if
\[ \lim_{n} |\ell_n(z)|^{1/n} = \exp\{g_{\Omega_0}(z;\infty)\}\,,
\]
uniformly on compact subsets of the complement of the smallest
interval containing the support, $\sop(\sigma_0),$ of $\sigma_0 $
and $g_{\Omega_0}(\cdot;\infty)$ denotes Green's function for the
region $\Omega_0 = \mathbb{C} \setminus \sop(\sigma_0)$ with
singularity at $\infty.$ Analogously, one defines regularity for
the other measures $\sigma_1,\ldots,\sigma_m\,.$ In the sequel, we
write $(\sigma_0;\sigma_1,\ldots,\sigma_m) \in \mbox{\bf Reg}$ to
mean that $\sigma_k \in \mbox{\bf Reg}, k=0,\ldots,m\,.$ The
system $(\sigma_1,\ldots,\sigma_m)$ will be used to construct the
Angelesco system of functions whereas $\sigma_0$ will determine
the system of orthogonal polynomials with respect to which the
Fourier expansions will be taken. Therefore, for all $0 \leq j,k
\leq m\,,$ we assume that
\[ \Delta_j \cap \Delta_k = \emptyset\,, \qquad j \neq k\,.
\]

In Theorems 1 and 2 below, we  find the rate of convergence of the
$|{\bf n}|$th root of the error of approximation of the functions
$\widehat{\sigma}_k$ by linear and non-linear Fourier-Pad\'e
approximants, respectively. The answers are given in terms of
extremal solutions of certain vector valued equilibrium problems
for the logarithmic potential. Before stating Theorems 1 and 2, we
need to introduce some notation and results from potential theory.

Let $F_k, k=1,\ldots,N,$ be (not necessarily distinct) closed
bounded intervals of the real line and $\mathcal{C} = (c_{j,k})$
be a real, positive definite, symmetric matrix of order $N$.
$\mathcal{C}$ will be called the interaction matrix. By
$\mathcal{M}_1(F_k), k =1,\ldots,N,$ we denote the subclass of
probability measures of $\mathcal{M}(F_k)$ and
\[\mathcal{M}_1= \mathcal{M}_1(F_1)
\times \cdots \times \mathcal{M}_1(F_{N})  \,.
\]
Given a vector measure $\mu \in \mathcal{M}_1$ and $j=1,\ldots,N,$
we define the combined potential
\begin{equation}
\label{combpot} W^{\mu}_j(x) = \sum_{k=1}^{N} c_{j,k} V^{\mu_k}(x)
\,,\qquad x \in \Delta_j\,,
\end{equation}
where
\[ V^{\mu_k}(x) = \int \log \frac{1}{|x-t|} \,d\mu_k(t)\,
\]
denotes the standard logarithmic potential of $\mu_k$. We denote
\[ w_j^{\mu} = \inf \{W_j^{\mu}(x): x \in \Delta_j\} \,.
\]

In Chapter 5 of \cite{niso} (see Propositions 4.5, 4.6, and
Theorem 4.1) the authors prove (we state the result in a form
convenient for our purpose)

\begin{lem} \label{niksor}
Let $\mathcal{C}$ be a real, positive definite, symmetric matrix
of order $N$. If there exists  $\overline{\mu} =
(\overline{\mu}_1,\ldots,\overline{\mu}_{N})\in \mathcal{M}_1$
such that for each $j=1,\ldots,N$
\[ W_j^{\overline{\mu}} (x) = w_j^{\overline{\mu}}\,, \qquad x \in \sop (\overline{\mu}_j)\,,\]
then $\overline{\mu}$ is unique. Moreover, if $c_{j,k} \geq 0$
when $F_j \cap F_j \neq \emptyset$ then $\overline{\mu}$ exists.
\end{lem}

The vector measure $\overline{\mu} \in \mathcal{M}_1$ is called
the equilibrium solution for the vector potential problem
determined by $\mathcal{C}$ on the system of intervals $F_j\,, j =
1,\ldots,N\,.$

In the sequel $\Lambda = \Lambda(p_1,\ldots,p_m) \subset
\mathbb{Z}_+^m$ is an infinite system of distinct multi-indices
such that
\[ \lim_{ {\bf n} \in
\Lambda} \frac{n_j}{|{\bf n}|} = p_j \in (0,1)\,, \qquad
j=1,\ldots,m\,.
\]

Let us define the block matrix
\begin{displaymath}
\mathcal{C}_1= \left( \begin{array}{cc}  \mathcal{C}_{1,1} &
\mathcal{C}_{1,2} \\ \mathcal{C}_{2,1} & \mathcal{C}_{2,2}
\end{array}\right)\,,
\end{displaymath}
where
\begin{displaymath}
\mathcal{C}_{1,1}= \left( \begin{array}{cccc}  2p_1^2 & p_1p_2 &
\cdots & p_1p_m \\ p_2p_1 & 2p_2^2 & \cdots & p_2p_m \\
\vdots & \vdots & \ddots & \vdots \\
p_mp_1 & p_mp_2 & \cdots & 2p_m^2
\end{array}\right)\,,\end{displaymath}
and $\mathcal{C}_{1,2}, \mathcal{C}_{2,1}, \mathcal{C}_{2,2} $ are
diagonal matrices given by
\begin{displaymath}
\mathcal{C}_{1,2}= \mathcal{C}_{2,1} =
\mbox{diag}\{-p_1(1+p_1),-p_2(1+p_2),\cdots,-p_m(1+p_m)\} \,,
\end{displaymath}
and
\begin{displaymath}
\mathcal{C}_{2,2} =
\mbox{diag}\{2(1+p_1)^2,2(1+p_2)^2,\cdots,2(1+p_m)^2\} \,.
\end{displaymath}
$\mathcal{C}_1$ satisfies all the assumptions of Lemma
\ref{niksor} on the system of intervals $F_j = \Delta_j,
j=1,\ldots,m, F_j = \Delta_0, j = m+1,\ldots,2m,$ including
$c_{j,k} \geq 0$ when $F_j \cap F_j \neq \emptyset.$ The only
non-trivial property is its positive definiteness and we shall
prove this in Section 2. Let
$\overline{\mu}=\overline{\mu}(\mathcal{C}_1) $ be the equilibrium
solution for the corresponding vector potential problem.  We have

\begin{thm}\label{teoprin}
Let $({\sigma}_0;\sigma_1,\ldots,{\sigma}_m) \in \mbox{\bf Reg}$
and consider the sequence of multi-indices $\Lambda =
\Lambda(p_1,\ldots,p_m)$.  Let $\left(\frac{P_{{\bf n},1}}{Q_{\bf
n}},\ldots,\frac{P_{{\bf n},m}}{Q_{\bf n}}\right), {\bf n} \in
\Lambda, $ be the associated sequence of linear Fourier-Pad\'e
approximants for the Angelesco system of functions
$(\widehat{\sigma}_1,\ldots,\widehat{\sigma}_m)$ with respect to
$\sigma_0$. Then,
\begin{equation}\label{conver1}
\lim_{{\bf n} \in
\Lambda}\left|\widehat{\sigma}_j(z)-\frac{P_{{\bf n},j}(z)}{Q_{\bf
n}(z)}\right|^{1/|{\bf n}|}=G_j(z)\,, \qquad j=1,\ldots,m\,,
\end{equation}
uniformly on each compact subset of $\overline{\C}\setminus
(\cup_{k=0}^m \Delta_k)$, where
\[G_j(z)=\exp\left((W_j^{\overline{\mu}}(z)- \omega_j^{\overline{\mu}})/p_j\right)\,,\]
$\overline{\mu} = \overline{\mu}(\mathcal{C}_1),$ and the combined
potentials $W_j^{\overline{\mu}}$ are defined by $(\ref{combpot})$
using $\mathcal{C}_1$.
\end{thm}

Set
\[
G_j^{\pm}=\{x\in\overline{\C}\setminus(\cup_{k=0}^{m}
\Delta_k):\,\pm\left( \omega_j^{\overline{\mu}}-
W_j^{\overline{\mu}}(x)\right)>0\}.
\]
An immediate consequence of Theorem \ref{teoprin} is

\begin{cor}\label{convergencia-lineal}
Under the assumptions of Theorem \ref{teoprin},
\[ \lim_{{\bf n} \in \Lambda} \frac{P_{{\bf
n},j}}{Q_{\bf n}}=\widehat{\sigma}_j \,, \qquad j=1,\ldots,m\,,\]
uniformly on compact subsets of $G_j^+$ and diverges to infinity
at each point of $G_j^-$.
\end{cor}

Non-linear Fourier Pad{\'e} approximants require the solution of a
different vector potential equilibrium problem. Let
\begin{displaymath}
\mathcal{C}_2= \left( \begin{array}{cc}  \mathcal{C}_{1,1} &
\mathcal{C}_{1,2} \\ \mathcal{C}_{2,1} & \mathcal{C}_{2,2}^2
\end{array}\right)\,,
\end{displaymath}
where $\mathcal{C}_{1,1}, \mathcal{C}_{1,2},\mathcal{C}_{2,1}$ are
as before and
\[ \mathcal{C}_{2,2}^2 = \left(
\begin{array}{cccc}
\frac{2m(1+p_1)^2}{m+1} & \frac{-2(1+p_1)(1+p_2)}{m+1} & \cdots &
\frac{-2(1+p_1)(1+p_m)}{m+1} \\
\frac{-2(1+p_2)(1+p_1)}{m+1} & \frac{2m(1+p_2)^2}{m+1} & \cdots &
\frac{-2(1+p_2)(1+p_m)}{m+1} \\
\vdots & \vdots & \ddots & \vdots \\
\frac{-2(1+p_m)(1+p_1)}{m+1} & \frac{-2(1+p_m)(1+p_2)}{m+1} &
\cdots & \frac{2m(1+p_m)^2}{m+1}
\end{array}
\right)\,.
\]
$\mathcal{C}_2$ is a real, positive definite, symmetric matrix of
order $2m$. We take the system of intervals $F_j = \Delta_j,
j=1,\ldots,m, F_j = \Delta_0, j = m+1,\ldots,2m.$ $\mathcal{C}_2$
does not satisfy that $c_{j,k} \geq 0$ when $F_j \cap F_j \neq
\emptyset.$ In Theorem \ref{asint2} of Section 3, we prove that
the corresponding equilibrium problem has at least one solution
and that $\mathcal{C}_2$ is positive definite. Therefore,
according to Lemma \ref{niksor} the solution is unique.   Let
$\overline{\mu}=\overline{\mu}(\mathcal{C}_2) $ be the equilibrium
solution for the corresponding vector potential problem.  In Lemma
\ref{existe} we show that for each ${\bf n} \in \mathbb{Z}_+^m$
there exists at least one non-linear Fourier- Pad\'e approximant
but we have not been able to prove that it is unique. We have

\begin{thm}\label{teoprin2}
Let $({\sigma}_0;\sigma_1,\ldots,{\sigma}_m) \in \mbox{\bf Reg}$
and consider the sequence of multi-indices $\Lambda =
\Lambda(p_1,\ldots,p_m)$.  Let $\left(\frac{S_{{\bf n},1}}{T_{\bf
n}},\ldots,\frac{S_{{\bf n},m}}{T_{\bf n}}\right), {\bf n} \in
\Lambda, $ be an associated sequence of non-linear Fourier-Pad\'e
approximants for the Angelesco system of functions
$(\widehat{\sigma}_1,\ldots,\widehat{\sigma}_m)$ with respect to
$\sigma_0$. Then,
\begin{equation}\label{conver2}
\lim_{ {\bf n} \in
\Lambda}\left|\widehat{\sigma}_j(z)-\frac{S_{{\bf n},j}(z)}{T_{\bf
n}(z)}\right|^{1/|{\bf n}|}=H_j(z)\,, \qquad j=1,\ldots,m\,,
\end{equation}
uniformly on each compact subset of $\overline{\C}\setminus
(\cup_{k=0}^m \Delta_j)$, where
\[ H_j(z) =  \exp\left((W_j^{\overline{\mu}}(z)-
\omega_j^{\overline{\mu}})/p_j\right)\,,
\]
$\overline{\mu} = \overline{\mu}(\mathcal{C}_2)$ and the combined
potentials $W_j^{\overline{\mu}}$ are defined by $(\ref{combpot})$
using $\mathcal{C}_2$.
\end{thm}

Notice that the limit only depends on $\Lambda$ and not on the
non-linear Fourier-Pad\'e approximants selected (in case that they
were not uniquely determined). Set
\[
H_j^{\pm}=\{x\in\overline{\C}\setminus(\cup_{j=0}^{m}
\Delta_j:\,\pm\left( \omega_j^{\overline{\mu}}-
W_j^{\overline{\mu}}(x)\right)>0\}.
\]
As a consequence of Theorem \ref{teoprin2}, we obtain

\begin{cor}\label{convergencia-lineal}
Under the assumptions of Theorem \ref{teoprin2},
\[ \lim_{{\bf n} \in \Lambda}
 \frac{S_{{\bf n},j}}{T_{\bf n}} = \widehat{\sigma}_j\,, \qquad j
 =1,\ldots,m\,,
 \]
uniformly on compact subsets of $H_j^+$ and diverges to infinity
at each point of $H_j^-$.
\end{cor}

Section \ref{teo-pot} is dedicated to the proof of Theorem
\ref{teoprin} and Section \ref{demostracion} to that of Theorem
\ref{teoprin2}. Section \ref{justi} is dedicated to the
justification of Lemma \ref{niksor} as stated here since in
\cite{niso} the assumption $c_{j,k} \geq 0$ if $F_j \cap F_k \neq
\emptyset$ is assumed in general. In the sequel, we maintain the
notation introduced above.

\section{ Proof of Theorem \ref{teoprin}}\label{teo-pot}

From the definition of the linear Fourier-Pad\'e approximant
immediately follows that for each $j=1,\ldots,m$
\begin{equation}\label{four-pade-lineal}
\int x^k(Q_{\bf n}(x)\widehat{\sigma}_j(z)-P_{{\bf
n},j}(x))d\sigma_0(x)=0,\quad k=0,\ldots, |{\bf n}|+n_j-1\,.
\end{equation}
Since the function  $Q_{\bf n}(z)\widehat{\sigma}_j(z)-P_{{\bf
n},j}(z)$ is continuous on $\Delta_0$, from
(\ref{four-pade-lineal}) we have that $Q_{\bf
n}(z)\widehat{\sigma}_j(z)-P_{{\bf n},j}(z)$ has at least $|{\bf
n}|+n_j$ sign changes on $\Delta_0$.

Let $W_{{\bf n},j}$ be the monic polynomial whose zeros are the
points where $Q_{\bf n}(z)\widehat{\sigma}_j(z)-P_{{\bf n},j}(z)$
changes sign on the interval $\Delta_0$. Obviously, $\deg W_{{\bf
n},j} \geq |{\bf n}|+n_j$ and
\[
\frac{Q_{\bf n}(z)\widehat{\sigma}_j(z)-P_{{\bf n},j}(z)}{W_{{\bf
n},j}(z)}\in \mathcal{H}(\overline{\C}\setminus \sop(\sigma_j)),
\quad j=1,\ldots,m\,,
\]
is analytic on the indicated region. Thus, linear Fourier-Pad\'e
approximants satisfy interpolation conditions on $\Delta_0$. A
similar statement holds for the non-linear Fourier-Pad\'e
approximants. In our proofs, we will use certain orthogonality
relations   satisfied by vector rational interpolants.

\begin{lem}\label{ortogonalidad} Let $(\widehat{\sigma}_1,\ldots,\widehat{\sigma}_m)$
be an Angelesco system, ${\bf n} = (n_1,\ldots,n_m) \in
\mathbb{Z}_+^m,$ and $(w_{{\bf n},1},\ldots,w_{{\bf n},m})$ a
system of polynomials such that $\deg w_{{\bf n},j} \geq |{\bf
n}|+n_j, j=1,\ldots,m$ whose zeros lie on an interval $\Delta_0,
\Delta_0 \cap \Delta_j = \emptyset, j= 1,\ldots,m $. Let
$(\frac{p_{{\bf n},1}}{q_{\bf n}},\ldots,\frac{p_{{\bf
n},m}}{q_{\bf n}})$ be a vector rational function such that $\deg
p_{{\bf n},j} \leq |{\bf n}|-1, j=1,\ldots,m, \deg q_{\bf n} \leq
|{\bf n}|, q_{\bf n} \not\equiv 0,$ and
\begin{equation}\label{analiticidad}
\frac{q_{\bf n}(z)\widehat{\sigma}_j(z)-p_{{\bf n},j}(z)}{w_{{\bf
n},j}(z)}\in \mathcal{H}(\overline{\C}\setminus \sop(\sigma_j)),
\quad j=1,\ldots,m\,.
\end{equation}
Then
\begin{equation}\label{orto-total}
\int x^k\frac{q_{\bf n}(x)}{w_{{\bf n},j}(x)}d\sigma_j(x)=0\,,
 \quad k=0,\,1,\ldots,n_j-1\,, \quad j = 1,\ldots, m\,.
\end{equation}
Consequently, $\deg q_{\bf n} = |{\bf n}|$  with exactly $n_j$
simple zeros in the interior of $\Delta_j$  (in connection with
intervals of the real line, the interior refers to the Euclidean
topology of the real line) and $\deg w_{{\bf n},j} = |{\bf
n}|+n_j, j=1,\ldots,m$. Let $q_{\bf n}= q_{{\bf
n},j}\widetilde{q}_{{\bf n},j},$ where $q_{{\bf n},j}$ is the
monic polynomial whose zeros are those of $q_{\bf n}$ lying in the
interior of $\Delta_j.$ Then
\begin{equation}\label{resto1}
\widehat{\sigma}_j(z)-\frac{p_{{\bf n},j}(z)}{q_{\bf
n}(z)}=\frac{w_{{\bf n},j}(z)}{q_{{\bf
n},j}^2(z)\widetilde{q}_{{\bf n},j}(z)} \int \frac{q_{{\bf
n},j}^2(x)}{z-x}\frac{\widetilde{q}_{{\bf n},j}(x)}{w_{{\bf
n},j}(x)}d\sigma_j(x)\,.
\end{equation}
\end{lem}

{\bf Proof.} Notice that (\ref{analiticidad}) and the assumption
on the degrees of the polynomials $q_{\bf n}, p_{{\bf n},j},$ and
$w_{{\bf n},j}$ imply that for $j=1,\ldots,m,$ and
$k=0,\ldots,n_j-1,$
\[
z^k\frac{q_{\bf n}(z)\widehat{\sigma}_j(z)-p_{{\bf
n},j}(z)}{w_{{\bf
n},j}(z)}=\mathcal{O}\left(\frac{1}{z^2}\right),\quad
z\to\infty\,.
\]
Let  $\Gamma_j$ be a closed, smooth, Jordan curve that surrounds
$\Delta_j$ such that all the intervals $\Delta_i, i\ne j, i =
0,\ldots,m$, lie in the unbounded connected component of the
complement of $\Gamma_j$. By Cauchy's Theorem, Cauchy's Integral
Formula and Fubini's Theorem, it follows that
\[
0=\frac{1}{2\pi i}\int_{\Gamma_j}z^k\frac{q_{\bf
n}(z)\widehat{\sigma}_j(z)-p_{{\bf n},j}(z)}{w_{{\bf n},j}(z)}d z
=
\]
\[
\frac{1}{2\pi i}\int_{\Gamma_j}z^k\frac{q_{\bf
n}(z)\widehat{\sigma}_j(z)}{w_{{\bf n},j}(z)}d z-\frac{1}{2\pi
i}\int_{\Gamma_j}z^k\frac{p_{{\bf n},j}(z)}{w_{{\bf n},j}(z)}d z =
\]
\[
\frac{1}{2\pi i}\int_{\Gamma_j}z^k\frac{q_{\bf n}(z)}{w_{{\bf
n},j}(z)}\int \frac{d\sigma_j(x)}{z-x}\,d z=\int x^k\frac{q_{\bf
n}(x)}{w_{{\bf n},j}(x)}d\sigma_j(x),
\]
for $k=0,\,1,\ldots,n_j-1$ and $j=1,\ldots,m$. Therefore,
(\ref{orto-total}) follows.

Using standard arguments of orthogonality, from (\ref{orto-total})
we obtain that $q_{\bf n}$ must have at least $n_j$ sign changes
in the interior of $\Delta_j$ and, consequently, at least $n_j$
zeros of odd multiplicity. Since $\deg q_{\bf n} \leq |{\bf n}|,$
we have that $\deg q_{\bf n} = |{\bf n}|$, that all its zeros are
simple and they are distributed in such a way that exactly $n_j$
lie in the interior of $\Delta_j$.

Assume that $\deg w_{{\bf n},j} > |{\bf n}| + n_j$ for some $j$.
Then
\[
z^k\frac{q_{\bf n}(z)\widehat{\sigma}_j(z)-p_{{\bf
n},j}(z)}{w_{{\bf
n},j}(z)}=\mathcal{O}\left(\frac{1}{z^2}\right),\quad
z\to\infty\,, \qquad k = 0,\ldots,n_j\,.
\]
This implies that (\ref{orto-total}) holds for all $k =
0,\ldots,n_j$. In turn, this means that $q_{\bf n}$ has at least
$n_j +1$ zeros in the interior of $\Delta_j$ against what was just
proved. Therefore, $\deg w_{{\bf n},j} = |{\bf n}| + n_j.$

Set $q_{\bf n}= q_{{\bf n},j}\widetilde{q}_{{\bf n},j},$ where
$q_{{\bf n},j}$ is the monic polynomial whose zeros are those of
$q_{\bf n}$ lying in the interior of $\Delta_j$. Notice that
$\widetilde{q}_{{\bf n},j}d\sigma_j/w_{{\bf n},j}$ is a real
measure with constant sign on $\Delta_j$. For future reference,
notice that with this notation the orthogonality relations
(\ref{orto-total}) may be expressed as
\begin{equation} \label{relortvar}
\int x^kq_{{\bf n},j}(x)|\widetilde{q}_{{\bf
n},j}(x)|\frac{d\sigma_j(x)}{|w_{{\bf n},j}(x)|}=0,\quad
k=0,\,1,\ldots,n_j-1\,.
\end{equation}
Hence, for each $j=1,\ldots,m,$ $q_{{\bf n},j}$ is the monic
orthogonal polynomial of degree $n_j$ with respect to the varying
measure $\frac{|\widetilde{q}_{{\bf n},j}|}{|w_{{\bf
n},j}|}d\sigma_j\,.$

Notice that
\[
 \frac{[q_{{\bf n},j}(q_{\bf n}\widehat{\sigma}_j-p_{{\bf
n},j})](z)}{w_{{\bf
n},j}(z)}=\mathcal{O}\left(\frac{1}{z}\right),\quad z\to\infty\,.
\]
Choose  $\Gamma_j$ as before. Using Cauchy's integral formula,
Cauchy's Theorem, and Fubini's Theorem, we obtain that for each
$j=1,\ldots,m$
\[\frac{[q_{{\bf n},j}(q_{\bf n}\widehat{\sigma}_j-p_{{\bf
n},j})](z)}{w_{{\bf n},j}(z)} = \frac{1}{2\pi i}\int_{\Gamma_j}
\frac{[q_{{\bf n},j}(q_{\bf n} \widehat{\sigma}_j -p_{{\bf
n},j})](\zeta)}{w_{{\bf n},j}(\zeta)(z - \zeta)} d \zeta =
\]
\[ \int \frac{1}{2\pi i}\int_{\Gamma_j}
\frac{(q_{{\bf n},j}q_{\bf n})(\zeta)}{w_{{\bf n},j}(\zeta)(z -
\zeta)(\zeta -x)} d \zeta d \sigma_j(x)= \int \frac{q_{{\bf
n},j}^2(x)}{z-x}\frac{\widetilde{q}_{{\bf n},j}(x)}{w_{{\bf
n},j}(x)}d\sigma_j(x)\,,
\]
which is equivalent to (\ref{resto1}). \hfill $\Box$

The vector rational function $(\frac{p_{{\bf n},1}}{q_{\bf
n}},\ldots,\frac{p_{{\bf n},m}}{q_{\bf n}})$ is called a
multipoint vector Pad\'e approximant of the Angelesco system
$(\widehat{\sigma}_1,\ldots,\widehat{\sigma}_m)$. According to
Lemma \ref{ortogonalidad} a necessary condition for their
existence is that $\deg w_{{\bf n},j} \leq |{\bf n}|+n_j,
j=1,\ldots,m$. Solving a homogeneous linear system of equations
one sees that this condition is also sufficient. When $\deg
w_{{\bf n},j} = |{\bf n}|+n_j, j=1,\ldots,m$ uniqueness follows
because then $\deg q_{\bf n} = |{\bf n}|$ as we have seen.

\begin{rem} Applying this Lemma to linear Fourier-Pad\'e
approximants, we have that $\deg(Q_{\bf n})=|{\bf n}|$. Thus, for
each ${\bf n } \in \mathbb{Z}_+^m$, they are uniquely determined
as claimed.
\end{rem}

Let us return to linear Fourier-Pad\'e approximants. In this case,
$Q_{{\bf n}} = q_{{\bf n}}, Q_{{\bf n},j} = q_{{\bf n},j},
\widetilde{Q}_{{\bf n},j} = \widetilde{q}_{{\bf n},j}$ and
$W_{{\bf n},j} = w_{{\bf n},j}$.

\begin{lem}\label{prop-wn}
 For each $ j=1, \ldots,m,$ and $k=0,\ldots,|{\bf n}|+n_j-1$
\begin{equation}\label{ortoW1}
\int  t^k\frac{W_{{\bf n},j}(t)}{|Q_{{\bf n},j}(t)|}\left(\int
\frac{Q_{{\bf n},j}^2(x)}{|t-x|}\frac{|\widetilde{Q}_{{\bf
n},j}(x)|}{|W_{{\bf n},j}(x)|}d\sigma_j(x)\right)d\sigma_0(t)=0\,.
\end{equation}
Moreover, $\deg W_{{\bf n},j} = |{\bf n}|+n_j, j=1, \ldots,m;$
that is, $Q_{\bf n}(z)\widehat{\sigma_j}(z)-P_{{\bf n},j}(z)$ has
exactly $|{\bf n}|+n_j$ sign changes in the interior of
$\Delta_0$.
\end{lem}

{\bf Proof.} From (\ref{resto1}) and the definition of the linear
Fourier-Pad\'e approximant, (\ref{ortoW1}) follows directly. The
assertion concerning the degree of $W_{{\bf n},j}$ is also
contained in Lemma \ref{ortogonalidad}. \hfill $\Box$

Let $\{\mu_l\} \subset \mathcal{M}(\mathcal{K})$ be a sequence of
measures, where $\mathcal{K}$ is a compact subset of the complex
plane and $\mu \in \mathcal{M}(\mathcal{K})$. We write
\[ *\lim_l \mu_l = \mu\,, \qquad \mu \in
\mathcal{M}(\mathcal{K})\,,\] if for every continuous function $f
\in \mathcal{C}( \mathcal{K})$
\[ \lim_l \int f d \mu_l
 = \int f d\mu\,;\]
that is, when the sequence of measures converges to $\mu$ in the
weak star topology. Given a polynomial $q_l$ of degree $l \geq 1$,
we denote the associated normalized zero counting measure by
\[ \nu_{q_l} = \frac{1}{l} \sum_{q_l(x) = 0} \delta_x \,,\]
where $\delta_x$ is the Dirac measure with mass $1$ at $x$ (in the
sum the zeros are repeated according to their multiplicity).

In order to prove our main results we need Theorem 3.3.3 of
\cite{stto}. We present it in the form stated in \cite{gora} which
is more adequate for our purpose. In \cite{gora}, it was proved
under stronger assumptions on the measure.

\begin{lem}\label{gonchar-rakhmanov}
Let $\{\phi_l\}, l \in \Gamma \subset \mathbb{Z}_+,$ be a sequence
of positive continuous functions on a bounded closed interval
$\Delta \subset\R , $ $\sigma \in {\bf Reg} \cap
\mathcal{M}(\Delta),$ and let $\{q_l\}, l \in \Gamma,$ be a
sequence of monic polynomials such that $\deg q_l = l$ and
\[
\int q_l(t)t^k\phi_l(t)d\sigma(t)=0,\quad k=0,\ldots, l-1.
\]
Assume that
\[
\lim_{l\in \Gamma}\frac{1}{2l}\log\frac{1}{|\phi_l(x)|}= v(x),
\]
uniformly on $\Delta$. Then
\[
*\lim_{l \in \Gamma}\nu_{q_l} = \overline{\nu},
\]
and
\[
\lim_{l\in\Gamma}\left(\int |q_l|^2\phi_l d\mu\right)^{1/{2l}}=
e^{-\omega},
\]
where $\overline{\nu} \in \mathcal{M}_1(\Delta_1)$ is the
equilibrium measure for the extremal problem
\[
V^{\overline{\nu}}(x)+v(x) \left\{ \begin{array}{l} = \omega,\quad
x
\in \sop(\overline{\nu}) \,, \\
\geq \omega, \quad x \in \Delta_1 \,,
\end{array} \right.
\]
in the presence of the external field $v$.
\end{lem}

Using this result, we can obtain the asymptotic limit distribution
of the zeros of the polynomials $Q_{{\bf n},j}$ and $W_{{\bf
n},j}$.

\begin{thm} \label{asint1} Let $({\sigma}_0;\sigma_1,\ldots,{\sigma}_m) \in \mbox{\bf
Reg}$ and consider the sequence of multi-indices $\Lambda =
\Lambda(p_1,\ldots,p_m)$. Then, for each $j=1,\ldots,m$
\[ *\lim_{{\bf n} \in \Lambda }{\nu}_{Q_{{\bf n},j}} =
\overline{\mu}_j\,, \quad\quad *\lim_{{\bf n} \in \Lambda
}\nu_{W_{{\bf n},j}} = \overline{\mu}_{m+j}\,,
\]
where $\overline{\mu} = \overline{\mu}(\mathcal{C}_1)  \in
\mathcal{M}_1$ is the vector equilibrium measures determined by
the matrix $ \mathcal{C}_1$ on the system of intervals $F_j =
\Delta_j, j= 1,\ldots,m\,,\,\, F_j = \Delta_0, j=
m+1,\ldots,2m\,.$
\end{thm}

{\bf Proof.} The unit ball in the cone of positive Borel measures
is weakly compact; therefore, it is sufficient to show that the
sequences of measures $\{\nu_{Q_{{\bf n},j}}\}$ and
$\{\nu_{W_{{\bf n},j}}\}, {\bf n} \in \Lambda,$ have only one
accumulation point which coincide with the components of the
vector measure $\overline{\mu}(\mathcal{C}_1) $ respectively. Let
$\Lambda^{\prime} \subset \Lambda$ be a subsequence of
multi-indices such that for each $j=1,\ldots,m$
\[
*\lim_{{\bf n} \in \Lambda^{\prime}}\nu_{Q_{{\bf n},j}} = \nu_j\,,
\quad\quad *\lim_{{\bf n} \in \Lambda^{\prime}}\nu_{W_{{\bf n},j}}
= \nu_{m+j}\,.
\]
(Notice that $\nu_j \in \mathcal{M}_1(\Delta_j), j=1,\ldots,m,$
and $\nu_j \in \mathcal{M}_1(\Delta_0), j=m+1,\ldots,2m.$)
Therefore,
\begin{equation} \label{pol-pot1}
\lim_{{\bf n} \in \Lambda^{\prime}}|Q_{{\bf
n},j}(z)|^{\frac{1}{n_j}}=\exp(-V^{\nu_j}(z)),
\end{equation}
uniformly on compact subsets of $\C\setminus \Delta_j $, and
\begin{equation} \label{pol-pot2}
\lim_{{\bf n} \in \Lambda^{\prime}} |W_{{\bf
n},j}(z)|^{\frac{1}{|{\bf n}|+n_j}}=\exp(-V^{\nu_{m+j}}(z)),
\end{equation}
uniformly on  compact subsets of $\C\setminus \Delta_0$.

For each fixed $j=1\ldots,m,$ the polynomials $Q_{{\bf n},j}$
satisfy the orthogonality relations (\ref{relortvar}). Using
(\ref{pol-pot1}) and (\ref{pol-pot2}) it follows that
\[
\lim_{{\bf n}\in \Lambda^{\prime}}\frac{1}{2n_j}\log\frac{|W_{{\bf
n},j}(x)|}{|\widetilde{Q}_{{\bf n},j}(x)|}=
-\frac{1+p_j}{2p_j}V^{\nu_{m+j}}(x)+ \sum_{k \neq j}
\frac{p_k}{2p_j}V^{\nu_k}(x),
\]
uniformly on $\Delta_j$. By Lemma \ref{gonchar-rakhmanov}, $\nu_j$
is the unique equilibrium measure for the extremal problem
\begin{equation}\label{prob-extre-11}  V^{\nu_j}(x)+ \sum_{k \neq j}
\frac{p_k}{2p_j}V^{\nu_k}(x) -\frac{1+p_j}{2p_j}V^{\nu_{m+j}}(x)
\geq \theta_{j},\quad x\in \Delta_j \,,
\end{equation}
with equality for all $x \in \sop(\nu_j)$.  Additionally,
\begin{equation}\label{asin-1}
\lim_{{\bf n}\in \Lambda^{\prime}} \left(\int |Q_{{\bf
n},j}(x)|^2\frac{|\widetilde{Q}_{{\bf
n},j}(x)|d\sigma_j(x)}{|W_{{\bf
n},j}(x)|}\right)^{\frac{1}{2n_j}}=e^{-\theta_j}.
\end{equation}

On the other hand, for each fixed $j=1,\ldots,m,$ the polynomials
$W_{{\bf n},j}$ satisfy the orthogonality relations (\ref{ortoW1})
and we can apply once more Lemma \ref{gonchar-rakhmanov}. Notice
that for all $t \in \Delta_0$
\begin{equation}\label{desrest1}
\int \frac{|Q_{{\bf n},j}^2(x)|}{|t-x|}\frac{|\widetilde{Q}_{{\bf
n},j}(x)|d\sigma_j(x)}{|W_{{\bf n},j}(x)|} \le\frac{1}{\delta_j}
\int |Q_{{\bf n},j}(x)|^2\frac{|\widetilde{Q}_{{\bf
n},j}(x)|d\sigma_j(x)}{|W_{{\bf n},j}(x)|},
\end{equation}
where  $\delta_j= \inf\{|t-x|: t \in \Delta_0,\, x \in \Delta_j\}$
and
\begin{equation}\label{desrest2}
\int \frac{|Q_{{\bf n},j}^2(x)|}{|t-x|}\frac{|\widetilde{Q}_{{\bf
n},j}(x)|d\sigma_j(x)}{|W_{{\bf n},j}(x)|}\ge \frac{1}{\delta_j^*}
\int |Q_{{\bf n},j}(x)|^2\frac{|\widetilde{Q}_{{\bf
n},j}(x)|d\sigma_j(x)}{|W_{{\bf n},j}(x)|},
\end{equation}
with $\delta_j^*=\max\{|t-x|:t \in \Delta_0,\,x\in \Delta_j\}$.
From (\ref{pol-pot1}), (\ref{pol-pot2}), (\ref{asin-1}),
(\ref{desrest1}), and (\ref{desrest2}), we obtain
\begin{multline*}
\lim_{{\bf n} \in \Lambda^{\prime}}\frac{1}{2(|{\bf
n}|+n_j)}\log\frac{|Q_{{\bf n},j}(x)|} {\int\frac{|Q_{{\bf
n},j}(t)|^2}{|x-t|}\frac{|\widetilde{Q}_{{\bf
n},j}(t)|d\sigma_j(t)}{|W_{{\bf n},j}(t)|}}
\\
=-\frac{p_j}{2(1+p_j)}V^{\nu_j}(x)+\frac{p_j}{1+p_j}\theta_j,
\end{multline*}
uniformly on $\Delta_0$.  Using Lemma \ref{gonchar-rakhmanov},
$\nu_{m+j}$ is the unique extremal solution for the equilibrium
problem
\begin{equation}\label{prob-extre-31}
V^{\nu_{m+j}}(x)-\frac{p_j}{2(1+p_j)}V^{\nu_j}(x)+\frac{p_j}{1+p_j}\theta_j
\geq  \theta_{m+j},\quad x\in \Delta_0\,,
\end{equation}
with equality for all $x \in \sop{(\nu_{m+j})}.$

Rewriting (\ref{prob-extre-11}) and (\ref{prob-extre-31})
conveniently, we see that the vector measure
$(\nu_1,\ldots,\nu_{2m}) \in \mathcal{M}_1$  is the unique
solution for the vector equilibrium problem determined by the
system of extremal problems
\begin{equation}\label{prob-extre-11*} 2p_j^2V^{\nu_j}(x)+ \sum_{k \neq
j}  p_jp_k V^{\nu_k}(x) -p_j(1+p_j) V^{\nu_{m+j}}(x) \geq
\omega_j,\quad x\in \Delta_j \,,
\end{equation}
($2p_j^2\theta_j = \omega_j$) with equality for all $x \in
\sop(\nu_j)$, and
\begin{equation}\label{prob-extre-31*}
2(1 + p_j)^2 V^{\nu_{m+j}}(x)- p_j(1+p_j) V^{\nu_j}(x) \geq
\omega_{m+j}\,,\quad x\in \Delta_0\,,
\end{equation}
with equality for all $x \in \sop{(\nu_{m+j})}.$ That is, it is
the equilibrium measure $\overline{\mu} \in \mathcal{M}_1$ for the
vector potential problem determined by $\mathcal{C}_1$ on the
system of intervals $F_j = \Delta_j, j = 1,\ldots,m,\, F_j =
\Delta_0, j=m+1,\ldots,2m.$ The condition $c_{j,k} \geq 0$ if $F_j
\cap F_k \neq \emptyset$ is fulfilled. According to Lemma
\ref{niksor}, this equilibrium vector measure is uniquely
determined if $\mathcal{C}_1$ is positive definite. Let us prove
this.

For $j \in \{1,\ldots,m\}$ the principle minor
$\mathcal{C}_1^{(j)}$ of order $j$ of $ \mathcal{C}_1$ is
\[ \mbox{det}(\mathcal{C}_1^{(j)}) = (p_1\cdots p_j)^2
\left| \begin{array}{cccc} 2 & 1 & \cdots & 1 \\
1 & 2 & \cdots & 1 \\
\vdots & \vdots & \ddots & \vdots \\
1 & 1 & \cdots & 2
\end{array}
 \right|_{j\times j} = (p_1\cdots p_j)^2(j+1) > 0 \,.
\]
For $j \in \{m+1,\ldots,2m\}$ the principle minor
$\mathcal{C}_1^{(j)}$ of order $j$ of $ \mathcal{C}_1$ can be
calculated as follows. For each $k=1,\ldots,m,$ factor out  $p_k$
from the $k$th row and $k$th column of $\mathcal{C}_1^{(j)}.$ From
the row and column $m+k, k=1,\ldots, j-m,$ factor out $1 + p_k$.
In the resulting determinant, for each $k=1,\ldots,j-m,$ add  the
$k$th row to the $(m+k)$th row and then to the resulting
determinant add the $k$th column to the $(m+k)$th column. We
obtain
\[ \mbox{det}(\mathcal{C}_1^{(j)}) = [p_1\cdots p_m(1+p_1)\cdots(1 + p_{j-m})]^2
\left| \begin{array}{cccc} 2 & 1 & \cdots & 1 \\
1 & 2 & \cdots & 1 \\
\vdots & \vdots & \ddots & \vdots \\
1 & 1 & \cdots & 2
\end{array}
 \right|_{m\times m} =
\]
\[ [p_1\cdots p_m(1+p_1)\cdots(1 + p_{j-m})]^2(m+1) > 0 \,.
\]
With this we conclude the proof. \hfill $\Box$

{\bf Proof of Theorem \ref{teoprin}.} From (\ref{resto1}),
(\ref{desrest1}), and (\ref{desrest2}), the asymptotic behavior of
the function $\widehat{\sigma}_j -\frac{P_{{\bf n},j}}{Q_{\bf n}}$
depends on the behavior of $W_{{\bf n},j}$, $Q_{{\bf n},j}$,
$\widetilde{Q}_{{\bf n},j}$, and $\gamma_{{\bf n},j}$, where
\begin{align*}
\frac{1}{\gamma_{{\bf n},j}^2}&=\min_Q \left\{\int
|Q(x)|^2\frac{|\widetilde{Q}_{{\bf n},j}(x)|d\sigma_j(x)}{|W_{{\bf
n},j}(x)|} :\,Q(x)=x^{n_j}+\cdots\right\}
\\
&=\int |Q_{{\bf n},j}(x)|^2\frac{|\widetilde{Q}_{{\bf
n},j}(x)|d\sigma_j(x)}{|W_{{\bf n},j}(x)|}.
\end{align*}

From Theorem \ref{asint1}, for each $j=1,\ldots,m,$ we have
\begin{equation}\label{vel-w1}
\lim_{{\bf n} \in \Lambda} |W_{{\bf n},j}(x)|^{1/{|{\bf
n}|}}=\exp\{-(1+p_j)V^{\overline{\mu}_{m+j}}(x)\},
\end{equation}
uniformly on compact subsets of $\C\setminus\Delta_0$, and
\begin{equation}\label{vel-q1}
\lim_{{\bf n} \in \Lambda} |Q_{{\bf n},j}(x)|^{1/{|{\bf
n}|}}=\exp\{-p_jV^{\overline{\mu}_j}(x)\},
\end{equation}
uniformly on compact subsets of $\C\setminus\Delta_j$, where
$\overline{\mu} = \overline{\mu}(\mathcal{C}_1)$ Using
(\ref{asin-1}) (see also parenthesis after
(\ref{prob-extre-11*})), it follows that
\begin{equation}\label{vel-res}
\lim_{|{\bf n}|\to\infty} \left(\frac{1}{\gamma_{{\bf
n},j}^2}\right)^{1/{|{\bf n}|}}=\exp \{ -2p_j\theta_j\} =
\exp\{-\omega_j/p_j\}\,.
\end{equation}

Combining (\ref{resto1}), (\ref{vel-w1}), (\ref{vel-q1}), and
(\ref{vel-res}), we conclude that (\ref{conver1}) holds true
uniformly on compact subsets of the indicated region. \hfill
$\Box$

\section{Proof of Theorem \ref{teoprin2}}\label{demostracion}

We begin by proving the existence of non-linear Fourier-Pad\'e
approximants.

\begin{lem} \label{existe} Given $(\sigma_0;\sigma_1,\ldots,\sigma_m),$
for each ${\bf n} \in \mathbb{Z}_+^m$ there exists an {\bf n}-th
non-linear Fourier-Pad\'e approximant of
$(\widehat{\sigma}_1,\ldots,\widehat{\sigma}_m)$ with respect to
$\sigma_0$.
\end{lem}

{\bf Proof.} In the proof  we make use of multipoint
Hermite-Pad\'e approximation.  Fix ${\bf n} \in \mathbb{Z}_+^m.$
For each $j \in \{1,\ldots,m\}$, choose an arbitrary  set of
$|{\bf n}|+n_j$ points contained in   $\Delta_0$
\[
X_{{\bf n},j}=(x_{{\bf n},j,1},\ldots,x_{{\bf n},j,|{\bf
n}|+n_j})\in \Delta_{{\bf n},j}\,,
\]
where
\[ \Delta_{{\bf n},j} = \{(x_{1},\ldots,x_{|{\bf
n}|+n_j}) \in \Delta_0^{|{\bf n}| + n_j}: x_1 \leq \cdots \leq
x_{|{\bf n}|+n_j} \}\,.
\]
Let
\[
w_{{\bf n},j}(x)=(x-x_{{\bf n},j,1})\cdots(x-x_{{\bf n},j,|{\bf
n}|+n_j})\,,
\]
and consider the simultaneous multipoint Pad\'e approximant which
interpolates the functions $\widehat{\sigma}_j, j=1,\ldots,m,$ at
the zeros of $w_{{\bf n},j}$ respectively. That is, $(p_{{\bf
n},1}/q_{\bf n},\ldots, p_{{\bf n},m}/q_{\bf n})$ is a vector
rational function such that $\deg (p_{{\bf n},j})\le |{\bf n}|-1$,
$j=1,\ldots,m$, $\deg(q_{\bf n})\le |{\bf n}|$, $q_{\bf
n}\not\equiv 0,$ and
\begin{equation}\label{def-mult}
\frac{q_{\bf n}\widehat{\sigma}_j-p_{{\bf n},j}}{w_{{\bf n},j}}\in
\mathcal{H}(\C\setminus\sop(\sigma_j)).
\end{equation}

From Lemma \ref{ortogonalidad} we have (\ref{orto-total}) and
(\ref{resto1}). Once we have determined $q_{\bf n}\,,$ for each
$j=1,\ldots,m,$ we define the monic polynomial $\Omega_{{\bf
n},j}\,, \deg(\Omega_{{\bf n},j}) = |{\bf n}|+n_j,$ by the
orthogonality relations
\begin{equation}\label{orto-aux1}
\int y^k\Omega_{{\bf n},j}(y)\left(\frac{1}{q_{{\bf
n},j}^2(y)\widetilde{q}_{{\bf n},j}(y)} \int \frac{q_{{\bf
n},j}^2(x)}{y-x}\frac{\widetilde{q}_{{\bf
n},j}(x)d\sigma_j(x)}{w_{{\bf n},j}(x)}\right)d\sigma_0(y)=0,
\end{equation}
$k=0,\ldots, |{\bf n}|+n_j-1$. For each $j=1,\ldots,m,$ these
relations determine a unique $\Omega_{{\bf n},j}$ since the
(varying) measures involved have constant sign on $\Delta_0\,.$

The polynomial $\Omega_{{\bf n},j}$ has exactly $|{\bf n}|+n_j$
 simple zeros in the interior of $\Delta_0\,.$ Set
\[
Y_{{\bf n},j}=(y_{{\bf n},j,1},\ldots,y_{{\bf n},j,|{\bf
n}|+n_j})\in \Delta_{{\bf n},j}\,,
\]
where $y_{{\bf n},j,1}<\cdots < y_{{\bf n},j,|{\bf n}|+n_j}$ are
the zeros of $\Omega_{{\bf n},j}$.

Since for each $j=1,\ldots,m,$ the distance between $\Delta_j$ and
$\Delta_0$ is greater than zero, the correspondence
\[
(X_{{\bf n},1},\ldots,X_{{\bf n},m}) \longrightarrow (Y_{{\bf
n},1},\ldots,Y_{{\bf n},m}),
\]
defines a continuous function  from $\Delta_{{\bf n},1} \times
\cdots \times \Delta_{{\bf n},m}$ into itself with the Euclidean
norm. The continuity of this function is an easy consequence of
the fact that $\Delta_0 \cap \Delta_j = \emptyset, j=1,\ldots,m.$
By Brouwer's fixed point Theorem (see page 364 of \cite{RS}) this
function has at least one fixed point. Choose a fixed point. Then,
$w_{{\bf n},j}=\Omega_{{\bf n},j}$, $j=1,\ldots,m\,.$ Consequently
(\ref{orto-aux1}) can be rewritten as
\begin{equation}\label{orto-aux2}
\int y^k w_{{\bf n},j}(y)\left(\frac{1}{q_{{\bf
n},j}^2(y)\widetilde{q}_{{\bf n},j}(y)} \int\frac{q_{{\bf
n},j}^2(x)}{y-x}\frac{\widetilde{q}_{{\bf
n},j}(x)d\sigma_j(x)}{w_{{\bf n},j}(x)}\right)d\sigma_0(y)=0\,,
\end{equation}
$k=0,\ldots, |{\bf n}|+n_j-1$, and taking into consideration
(\ref{resto1}) we obtain that for each $j=1,\ldots,m,$
\[
\int \left(\widehat{\sigma}_j(x)-\frac{p_{{\bf n},j}(x)}{q_{\bf
n}(x)}\right)x^kd\sigma_j(x)=0\,,\qquad k=0,\ldots,|{\bf n}|+n_j-1
\,.
\]
From the definition, it follows that   $(p_{{\bf n},1}/q_{\bf
n},\ldots,p_{n_m}/q_{\bf n})$ is an ${\bf n}$th non linear
Fourier-Pad\'e approximant for the Angelesco system, taking
$S_{{\bf n},j}=p_{{\bf n},j}$, $j=1,\ldots,m$, and $T_{\bf
n}=q_{\bf n}$.  \hfill $\Box$

Let $\left(\frac{S_{{\bf n},1}}{T_{\bf n}},\ldots,\frac{S_{{\bf
n},m}}{T_{\bf n}}\right)$ be any non-linear Fourier-Pad\'e
approximant with respect to the Angelesco system
$(\widehat{\sigma}_1,\ldots,\widehat{\sigma}_m)$. From ii') it
follows that  $ \widehat{\sigma}_j(z)-\frac{S_{{\bf
n},j}(z)}{T_{\bf n}(z)}$ has at least $|{\bf n}|+n_j$ sign changes
on $\Delta_0$. Let $W_{{\bf n},j}$ be the monic polynomial whose
zeros are the points where this function changes sign on
$\Delta_0$. Obviously, $\deg W_{{\bf n},j} \geq |{\bf n}|+n_j$ and
\begin{equation}\label{analiticidad1}
\frac{T_{\bf n}(z)\widehat{\sigma}_j(z)-S_{{\bf n},j}(z)}{W_{{\bf
n},j}(z)}\in \mathcal{H}(\overline{\C}\setminus \sop(\sigma_j)),
\quad j=1,\ldots,m\,,
\end{equation}
is analytic on the indicated region.  (These polynomials $W_{{\bf
n},j}$ do not coincide with those of the linear case.) Using Lemma
\ref{ortogonalidad} it follows that
\begin{equation}\label{x1}
\int x^k\frac{|T_{\bf n}(x)|}{|W_{{\bf n},j}(x)|}d\sigma_j(x)=0\,,
 \quad k=0,\,1,\ldots,n_j-1\,, \quad j = 1,\ldots, m\,.
\end{equation}
and
\begin{equation}\label{x2}
\widehat{\sigma}_j(z)-\frac{S_{{\bf n},j}(z)}{T_{\bf
n}(z)}=\frac{W_{{\bf n},j}(z)}{T_{{\bf
n},j}^2(z)\widetilde{T}_{{\bf n},j}(z)} \int \frac{T_{{\bf
n},j}^2(x)}{z-x}\frac{\widetilde{T}_{{\bf n},j}(x)}{W_{{\bf
n},j}(x)}d\sigma_j(x)\,,
\end{equation}
where $T_{{\bf n},j}$ is the monic polynomial whose zeros are the
$n_j$ zeros of $T_{\bf n}$ lying in the interior of $\Delta_j$.
Combining (\ref{x2}) with ii') we obtain
\begin{equation}\label{x3}
\int y^k W_{{\bf n},j}(y)\left(\frac{1}{T_{{\bf
n},j}^2(y)|\widetilde{T}_{{\bf n},j}(y)|} \int\frac{T_{{\bf
n},j}^2(x)}{|y-x|}\frac{|\widetilde{T}_{{\bf
n},j}(x)|d\sigma_j(x)}{|W_{{\bf n},j}(x)|}\right)d\sigma_0(y)=0\,.
\end{equation}

The proof of Theorem \ref{teoprin2} is similar to that of Theorem
\ref{teoprin}.  First, we study the asymptotic zero distribution
of the polynomials $T_{{\bf n},j}$ and $W_{{\bf n},j}$. Then, we
use this result to obtain the asymptotic behavior of the remainder
in the approximation.

\begin{thm} \label{asint2} Let $({\sigma}_0;\sigma_1,\ldots,{\sigma}_m) \in \mbox{\bf
Reg}$ and consider the sequence of multi-indices $\Lambda =
\Lambda(p_1,\ldots,p_m)$. Then, there exists a vector measure
$\overline{\mu} = (\overline{\mu}_1,\ldots,\overline{\mu}_{2m})
\in \mathcal{M}_1$ such that for each $j=1,\ldots,m$
\[ *\lim_{{\bf n} \in \Lambda }{\nu}_{T_{{\bf n},j}} = \overline{\mu}_j\,,
\quad\quad *\lim_{{\bf n} \in \Lambda }\nu_{W_{{\bf n},j}} =
\overline{\mu}_{m+j}\,.
\]
Moreover, $\overline{\mu} = \overline{\mu}(\mathcal{C}_2) $ is the
vector equilibrium measure determined by the matrix $
\mathcal{C}_2$ on the system of intervals $F_j = \Delta_j, j=
1,\ldots,m\,,\,\, F_j = \Delta_0, j= m+1,\ldots,2m\,.$
\end{thm}

{\bf Proof.} Let us show that the sequences of measures
$\{\nu_{T_{{\bf n},j}}\}$ and $\{\nu_{w_{{\bf n},j}}\}, {\bf n}
\in \Lambda,$ have only one accumulation point.  Let
$\Lambda^{\prime} \subset \Lambda$ be a subsequence of indices
such that for each $j=1,\ldots,m$
\[
*\lim_{{\bf n} \in \Lambda^{\prime}}\nu_{T_{{\bf n},j}} = \nu_j\,,
\quad\quad *\lim_{{\bf n} \in \Lambda^{\prime}}\nu_{W_{{\bf n},j}}
= \nu_{m+j}\,.
\]
(Notice that $\nu_j \in \mathcal{M}_1(\Delta_j), j=1,\ldots,m,$
and $\nu_j \in \mathcal{M}_1(\Delta_0), j=m+1,\ldots,2m.$)
Therefore,
\begin{equation} \label{pol-pot-no1}
\lim_{{\bf n} \in \Lambda^{\prime}}|T_{{\bf
n},j}(z)|^{\frac{1}{n_j}}=\exp(-V^{\nu_j}(z))\,,
\end{equation}
uniformly on compact subsets of $\C\setminus \Delta_j,$ and
\begin{equation} \label{pol-pot-no2}
\lim_{{\bf n} \in \Lambda^{\prime}}|W_{{\bf
n},j}(z)|^{\frac{1}{|{\bf n}|+n_j}}=\exp(-V^{\nu_{m+j}}(z)),
\end{equation} uniformly on compact subsets of $\C\setminus \Delta_0$.

As we have seen, $T_{{\bf n},j}$ is orthogonal with respect to the
varying measure $\frac{|\widetilde{T}_{{\bf n},j}|}{|W_{{\bf
n},j}|}d\sigma_j$. Using (\ref{pol-pot-no1}) and
(\ref{pol-pot-no2}), we obtain
\[
\lim_{{\bf n} \in\Lambda^{\prime}}\frac{1}{2n_j}\log\frac{|W_{{\bf
n},j}(x)|}{|\widetilde{T}_{{\bf n},j}(x)|}=
-\frac{1+p_j}{2p_j}V^{\nu_{m+j}}(x) + \sum_{k \neq j}
\frac{p_k}{2p_j}V^{\nu_k}(x)\,,
\]
uniformly in $\Delta_j$. By (\ref{x1}) and Lemma
\ref{gonchar-rakhmanov}, $\nu_j$ is the unique equilibrium measure
for the extremal problem
\begin{equation}\label{prob-11-no}
V^{\nu_j}(x)+\sum_{k \neq j}
\frac{p_k}{2p_j}V^{\nu_k}(x)-\frac{1+p_j}{2p_j}V^{\nu_{m+j}}(x)
\geq \eta_j\,,\qquad x\in \Delta_j \,,
\end{equation}
with equality for all $x \in \sop(\nu_j),$ and
\begin{equation}\label{asin-1-no}
\lim_{{\bf n} \in\Lambda^{\prime}} \left(\int |T_{{\bf n},j}^2(x)|
\frac{|\widetilde{T}_{{\bf n},j}(x)|d\sigma_j(x)}{|W_{{\bf
n},j}(x)|}\right)^{\frac{1}{2n_j}}=e^{-\eta_j}.
\end{equation}
These relations are completely similar to those obtained for the
linear case (see (\ref{prob-extre-11}) and (\ref{asin-1})). On the
other hand, $W_{{\bf n},j}$ satisfies the orthogonality relations
(\ref{x3}). We can apply once more Lemma \ref{gonchar-rakhmanov}
obtaining that, for each $j=1,\ldots,m,$ $\nu_{m+j}$ is the unique
equilibrium measure for the extremal problem
\begin{equation}\label{prob13-no}
V^{\nu_{m+j}}(x)- \frac{p_j}{1+p_j}V^{\nu_j}(x)- \sum_{k\neq j}
\frac{p_k}{2(1+p_j)}V^{\nu_k}(x) \geq  \eta_{m+j}\,,\qquad x\in
\Delta_{0} \,,
\end{equation}
with equality for all $x \in \sop(\nu_{m+j}).$ These relations
differ from those obtained for the linear case (see
(\ref{prob-extre-31}))

If we look at the matrix corresponding to this system of equations
we see that it is not symmetric. Let us rewrite the system as
follows. Multiply equations (\ref{prob-11-no}) times $2p_j^2$ and
we obtain for each $j=1,\ldots,m,$
\begin{equation} \label{c}
2p_j^2V^{\nu_j}(x)+\sum_{k \neq j} p_jp_k V^{\nu_k}(x)-p_j(1+p_j)
V^{\nu_{m+j}}(x) \geq 2\eta_jp_j^2 = w_j \,,\qquad x\in \Delta_j
\,,
\end{equation}
with equality for all $x \in \sop (\nu_j)$. With equations
(\ref{prob13-no}) we have to work harder. First, let us multiply
them times $2(1+p_j)$ thus obtaining for each $j=1,\ldots,m,$
\begin{equation} \label{d}
- 2p_j V^{\nu_j}(x)- \sum_{k\neq j}
 p_k V^{\nu_k}(x) + 2(1+p_j) V^{\nu_{m+j}}(x)\geq
\end{equation}
\[
2\eta_{m+j}(1+p_j) = \eta_{m+j}^{\prime}\,,\quad x\in \Delta_{0}
\,,
\]
with equality for all $x \in \sop (\nu_{m+j})$.

Let us show that in this second group of equations we have
equality for all $x\in \Delta_{0}$. In fact, notice that
$2(1+p_j)\nu_{m+j}$ is a measure on $\Delta_0$ of total mass equal
to $2(1+p_j)$. On the other hand
\[ 2p_j \nu_j + \sum_{k \neq j} p_k \nu_k
\]
is a measure of total mass $p_j + \sum_{k=1}^m p_k = 1 + p_j <
2(1+p_j)$ supported on the set $\cup_{k=1}^m \Delta_k$ which is
disjoint from $\Delta_0 .$ Therefore,
\[ 2(1+p_j)\nu_{m+j} =  ( 2p_j \nu_j + \sum_{k \neq j} p_k \nu_k
 )^{\prime}+ (1+p_j) \omega_{\Delta_0}\,,
\]
where $(\cdot)^{\prime}$ denotes the balayage onto $\Delta_0$ of
the indicated measure and $\omega_{\Delta_0}$ is the equilibrium
measure on $\Delta_0$ (without external field). since these two
measures are supported on all $\Delta_0$ so is their sum. Thus,
$\sop (\nu_{m+j}) = \Delta_0\,.$

The idea now is to take row transformations on the system of
equations (\ref{d}) to transform it conveniently. The matrix of
this system of equations is
\[\left(
\begin{array}{cccccccc}
-2p_1 & -p_2 & \cdots & -p_m & 2(1+p_1) & 0 & \cdots & 0 \\
-p_1 & -2p_2 & \cdots & -p_m & 0 & 2(1+p_2) & \cdots & 0 \\
\vdots & \vdots & \ddots & \vdots & \vdots & \vdots& \ddots & \vdots \\
-p_1 & -p_2 & \cdots & -2p_m & 0 & 0 & \cdots & 2(1+p_m)
\end{array} \right)\,.
\]
Since each column has a common factor we will carry out the
operations without the common factor and afterwards place them
back.  Thus in columns $k=1,\ldots,m$ we factor out $-p_k$ and in
columns $k=m+1,\ldots,2m$ we factor out $2(1+p_{k-m}),$
respectively. The resulting matrix is
\[\left(
\begin{array}{cccccccc}
2 & 1 & \cdots & 1 & 1 & 0 & \cdots & 0 \\
1 & 2  & \cdots & 1 & 0 & 1 & \cdots & 0 \\
\vdots & \vdots & \ddots & \vdots & \vdots & \vdots& \ddots & \vdots \\
1 & 1 & \cdots & 2 & 0 & 0 & \cdots & 1
\end{array} \right) = \left( \begin{array}{cc} \mathcal{B} & \mathcal{I} \end{array}
\right)\,,
\]
where $\mathcal{I}$ denotes the identity matrix of order $m$. We
know that the submatrix $\mathcal{B}$ is positive definite and
through row operations can be reduced to the identity. This is the
same as multiplying  $\left( \begin{array}{cc} \mathcal{B} &
\mathcal{I} \end{array} \right)$ on the left by
$\mathcal{B}^{-1}$. Doing this we obtain the block matrix
\[ \left( \begin{array}{cc}  \mathcal{I} & \mathcal{B}^{-1} \end{array}
\right)\,.
\]
It is easy to check that
\[  \mathcal{B}^{-1} =
\frac{1}{m+1}\left(
\begin{array}{cccc}
m & -1 & \ldots & -1 \\
-1 & m  & \ldots & -1 \\
\vdots & \vdots & \ddots & \vdots \\
-1 & -1 & \cdots & m
\end{array} \right)\,.
\]
Multiplying back the factors we extracted we obtain the matrix
\[\left(
\begin{array}{cccccccc}
-p_1 & 0 & \cdots & 0 & \frac{2m(1+p_1)}{m+1} & \frac{-2(1+p_2)}{m+1} & \cdots & \frac{-2(1+p_m)}{m+1} \\
0 & -p_2 & \cdots & 0 & \frac{-2(1+p_1)}{m+1} & \frac{2m(1+p_2)}{m+1} & \cdots & \frac{-2(1+p_m)}{m+1} \\
\vdots & \vdots & \ddots & \vdots & \vdots & \vdots& \ddots & \vdots \\
0 & 0 & \cdots & -p_m & \frac{-2(1+p_1)}{m+1} &
\frac{-2(1+p_2)}{m+1} & \cdots & \frac{2m(1+p_m)}{m+1}
\end{array} \right)\,.
\]
Therefore, the system of equations (\ref{d}) is equivalent to
\begin{equation} \label{e}
-p_j V^{\nu_j}(x) + \frac{2m(1+p_j)}{m+1} V^{\nu_{m+j}}(x) -
\end{equation}
\[
\sum_{k \neq j} \frac{2(1+p_k)}{m+1} V^{\nu_{m+j}}(x) =
\eta_{m+j}^{''}\,, \quad x \in \Delta_0\,,
\]
where
\[(\eta_{m+1}^{''},\ldots,\eta_{2m}^{''})^t=
\mathcal{B}^{-1}(\eta_{m+1}^{\prime},\ldots,\eta_{2m}^{\prime})^t.
\]
Finally, multiply the $j$th equation in (\ref{e}) times $(1+p_j)$
to obtain
\begin{equation} \label{f}
-p_j(1+p_j) V^{\nu_j}(x) + \frac{2m(1+p_j)^2}{m+1}
V^{\nu_{m+j}}(x) -
\end{equation}
\[\sum_{k \neq j} \frac{2(1+p_k)(1+p_j)}{m+1} V^{\nu_{m+j}}(x) =
\eta_{m+j}^{''}(1+p_j) = w_{m+j}\,, \quad x \in \Delta_0\,.
\]
The system of equilibrium problems defined by (\ref{c}) and
(\ref{f}) has the interaction matrix
\[ \mathcal{C}_2 = \left(
\begin{array}{cc}
\mathcal{C}_{1,1} & \mathcal{C}_{1,2} \\
 \mathcal{C}_{2,1} &
\mathcal{C}_{2,2}^2
\end{array}
\right)
\]
defined in Section 1. Thus, the corresponding equilibrium problem
has at least one solution given by $(\nu_1,\ldots,\nu_m)$.
According to Lemma 1, $(\nu_1,\ldots,\nu_m)$ is uniquely
determined if we prove that $\mathcal{C}_2$ is positive definite.

Let us show that $\mathcal{C}_2$ is positive definite. The first
$m$ principal minors of $\mathcal{C}_1$ and $\mathcal{C}_2$
coincide and we already know that they are positive. Let
$\mathcal{C}_2^{(j)}$ denote the principal minor of
$\mathcal{C}_2$ of order $j$ where $j \in \{ m+1,\ldots,2m\}$. For
each $k=1,\ldots,m,$ factor out $p_k$ from the $k$th row and $k$th
column of $\mathcal{C}_2^{(j)}.$ From the row and column $m+k,
k=1,\ldots, j-m,$ factor out $1 + p_k$. In the resulting
determinant, for each $k=1,\ldots,j-m,$ add  the $k$th row to the
$(m+k)$th row and then to the resulting determinant add the $k$th
column to the $(m+k)$th column. We obtain
\[ \mbox{det}(\mathcal{C}_2^{(j)}) = [p_1\cdots p_m(1+p_1)\cdots(1 +
p_{j-m})]^2 \times
\]
\[
\left| \begin{array}{cccccccc} 2 & 1 & \cdots & 1 & 1 & 1 & \cdots & 1 \\
1 & 2 & \cdots & 1  & 1 & 1 & \cdots & 1 \\
\vdots & \vdots & \ddots & \vdots & \vdots & \vdots & \ddots & \vdots \\
1 & 1 & \cdots & 2 & 1 & 1 & \cdots & 1 \\
1 & 1 & \cdots & 1 & \frac{2m}{m+1}  & \frac{m-1}{m+1} & \cdots & \frac{m-1}{m+1}\\
1 & 1 & \cdots & 1 & \frac{m-1}{m+1} & \frac{2m}{m+1} & \cdots & \frac{m-1}{m+1}\\
\vdots & \vdots & \ddots & \vdots & \vdots & \vdots & \ddots & \vdots \\
1 & 1 & \cdots & 1 & \frac{m-1}{m+1} & \frac{m-1}{m+1} & \vdots &
\frac{2m}{m+1}
\end{array}
 \right|
\]
In the determinant above, delete the row $m+1$ from the following
ones and in the resulting  determinant add to the column $m+1$
those after it and we get
\[
\left| \begin{array}{cccccccc} 2 & 1 & \cdots & 1 & j-m & 1 & \cdots & 1 \\
1 & 2 & \cdots & 1  & j-m & 1 & \cdots & 1 \\
\vdots & \vdots & \ddots & \vdots & \vdots & \vdots & \ddots & \vdots \\
1 & 1 & \cdots & 2 & j-m & 1 & \cdots & 1 \\
1 & 1 & \cdots & 1 &  \frac{(m+1)+(j-m)(m-1)}{m+1}  & \frac{m-1}{m+1} & \cdots & \frac{m-1}{m+1}\\
0 & 0 & \cdots & 0 & 0 & 1 & \cdots & 0 \\
\vdots & \vdots & \ddots & \vdots & \vdots & \vdots & \ddots & \vdots \\
0 & 0 & \cdots & 0 & 0 & 0 & \cdots & 1
\end{array}
 \right| =
\]
\[
\left| \begin{array}{ccccc} 2 & 1 & \cdots & 1 & j-m  \\
1 & 2 & \cdots & 1  & j-m  \\
\vdots & \vdots & \ddots & \vdots & \vdots   \\
1 & 1 & \cdots & 2 & j-m  \\
1 & 1 & \cdots & 1 &  j-m
\end{array}
 \right| +
\left| \begin{array}{ccccc} 2 & 1 & \cdots & 1 & 0  \\
1 & 2 & \cdots & 1  & 0 \\
\vdots & \vdots & \ddots & \vdots & \vdots   \\
1 & 1 & \cdots & 2 & 0  \\
1 & 1 & \cdots & 1 &  \frac{(m+1)-2(j-m)}{m+1}
\end{array}
 \right| =
\]
\[ (j-m) + (m+1) - 2(j-m) = 2m+1 -j > 0\,.
\]
With this we conclude the proof.  \hfill $\Box$

We are ready to prove Theorem \ref{teoprin2}.

{\bf Proof of Theorem \ref{teoprin2}.} From (\ref{x2}) the
asymptotic behavior of $\widehat{\sigma}_j(z)-\frac{S_{{\bf
n},j}(z)}{T_{\bf n}(z)}$ can be expressed in terms of that of the
sequences of polynomials $W_{{\bf n},j}$, $T_{{\bf n},j}$, and
$\zeta_{{\bf n},j},$ where
\[
\frac{1}{\zeta_{{\bf n},j}^2}=\min\left\{\int
|Q(x)|^2\frac{|\widetilde{T}_{{\bf n},j}(x)|d\sigma_j(x)}{|W_{{\bf
n},j}(x)|} :\,Q(x)=x^{n_j}+\cdots\right\}
\]
\[
=\int |T_{{\bf n},j}(x)|^2\frac{|\widetilde{T}_{{\bf
n},j}(x)|d\sigma_j(x)}{|W_{{\bf n},j}(x)|}\,.
\]
On account of Theorem \ref{asint2},   we have
\begin{equation}\label{vel-w1no}
\lim_{{\bf n} \in\Lambda} |W_{{\bf n},j}(z)|^{1/{|{\bf
n}|}}=\exp\{-(1+p_j)V^{\overline{\mu}_{m+j}}(z)\}\,,
\end{equation}
uniformly on compact subsets of $\C\setminus \Delta_0$, and
\begin{equation}\label{vel-t1no}
\lim_{{\bf n} \in\Lambda} |T_{{\bf n},j}^2(z)|^{1/{|{\bf
n}|}}=\exp\{-2p_jV^{\overline{\mu}_j}(z)\}\,,
\end{equation}
uniformly on compact subsets of $\C\setminus \Delta_j\,,$ where
$\overline{\mu} = \overline{\mu}(\mathcal{C}_2)$. Using
(\ref{asin-1-no}), we have
\begin{equation}\label{vel-res-no1}
\lim_{{\bf n} \in\Lambda} \left(\frac{1}{\zeta_{{\bf
n},j}^2}\right)^{1/{|{\bf n}|}}=\exp\{-2p_j\eta_j\} = \exp\{-
w_j/p_j\} \,.
\end{equation}
Combining (\ref{x2}), (\ref{vel-w1no}), (\ref{vel-t1no}), and
(\ref{vel-res-no1}), we obtain that (\ref{conver2}) holds true
uniformly on compact subsets of the indicated region. \hfill
$\Box$

\section{Comments on Lemma \ref{niksor}} \label{justi}

Let $\mathcal{M}(F_k), k =1,\ldots,N,$ be the class of all
measures on $\mathcal{M}(F_k)$ and
\[\mathcal{M}= \mathcal{M}_(F_1)
\times \cdots \times \mathcal{M}(F_{N})  \,.
\]
Define the mutual energy of two vector measures  $\mu^1,\mu^2 \in
{\mathcal{M}}$ by
\begin{equation}
\label{mutualenergy} J(\mu^1, \mu^2)= \sum_{j,k=1}^N \int \int
 c_{j,k} \ln \frac{1}{|z-x|}   d \mu^1_{j}(z) d
\mu^2_{k}(x).
\end{equation}
The energy of the vector measure $\mu \in {\mathcal{M}}$ is
\begin{equation}
\label{energy} J(\mu)= \sum_{j,k=1}^N c_{j,k}I(\mu_{j},\mu_{k})
\,,
\end{equation}
where
\begin{equation} \label{mutualenergy*}
I(\mu_{j},\mu_{k})=\int \int \ln \frac{1}{|z-x|} d \mu_{j}(z) d
\mu_{k}(x)\,.
\end{equation}
Given  a real symmetric positive definite matrix $\mathcal{C}$ and
$\mu \in \mathcal{M}$ define the combined potentials
$W^{\mu}_j,j=1,\ldots,N,$ as in the introduction, and the vector
potential $W^{\mu} = (W^{\mu}_1,\ldots,W^{\mu}_N).$ These formulas
may be rewritten as
\begin{equation}
\label{energiamutuapotencial} J(\mu^1, \mu^2)=\int W^{\mu^2}(z) d
\mu^1(z) \,,
\end{equation}
where
\[ \int W^{\mu^2}(z) d\mu^1(z) = \sum_{i=1}^m \int W_i^{\mu^2}(z)
d\mu_i^1(z)\,, \] and
\begin{equation}
\label{energiapotencial} J(\mu)=\int W^{\mu}(z) d \mu(z)\,.
\end{equation}
If $\mu, \mu^1, \mu^2 \in \mathcal{E}$ are vector charges whose
components have finite energy, the energy of a charge $J(\mu)$ and
the mutual energies of the charges $J(\mu^1,\mu^2)$ can be defined
analogously by formulas (\ref{energy}) and (\ref{mutualenergy}),
respectively.

If $c_{j,k} \geq 0$ when $F_j \cap F_k \neq \emptyset, j,k \in
\{1,\ldots,N\}$ the functionals $J(\mu^1, \mu^2)$ and $J(\mu)$ are
lower semicontinuous in the weak topology $\mathcal{M}$ (see
Proposition 5.4.1 in \cite{niso}). Consequently, the functional
$J(\mu)$ attains its minimum in $\mathcal{M}_1$. In Proposition
5.4.2, using a unitary decomposition of $\mathcal{C}$, the authors
prove that $J(\mu)$ is a nonsingular positive definite quadratic
form on the linear space $\mathcal{E}$. In Proposition 5.4.2, the
extra condition on the coefficients of $\mathcal{C}$ is not
needed. Therefore, we can say that if there is a minimizing vector
measure then it is unique.

Let $0 \leq \epsilon \leq 1$ and $\mu^1,\mu^2 \in
{\mathcal{M}_1}.$ Assume that the components of $\mu^1,\mu^2$ have
finite energy. Set $\widetilde{\mu}= \epsilon \mu^2 + (1-
\epsilon) {\mu}^1 \in {\mathcal{M}}_1$. It is  algebraically
straightforward to verify that
\begin{equation}
\label{epsilonal2} J(\widetilde{\mu}) - J(\mu^1) = \epsilon^2
J(\mu^2-\mu^1)+2\epsilon \int W^{ {\mu}^1}(x) d(\mu^2 -
{\mu}^1)(x) \,.
\end{equation}
Assume that $J(\mu^1)$ minimizes the energy functional. Dividing
by $\epsilon$ and letting $\epsilon$ tend to zero, it follows that
\begin{equation} \label{minimo}
 \int W^{ \mu^1}(x) d(\mu^2 - \mu^1)(x)  \geq 0.
\end{equation}
for all $\mu^2 \in \mathcal{M}_1$. Reciprocally, assume that
(\ref{minimo}) takes place for all $\mu^2 \in \mathcal{M}_1$, then
using (\ref{epsilonal2}) it follows that $\mu^1$ minimizes the
energy functional since $J(\mu^2-\mu^1) \geq 0$ for all
$\mu^1,\mu^2 \in \mathcal{C}$.

Now, let $\overline{\mu} \in \mathcal{M}_1$ be a solution of the
equilibrium potential problem determined by $\mathcal{C}$ on the
system of intervals $F_j\,, j = 1,\ldots,N\,.$ That is
\[ W_j^{\overline{\mu}} (x) = w_j^{\overline{\mu}}\,, \qquad x \in
\sop (\overline{\mu}_j)\,,
\]
where $w_j^{\overline{\mu}} = \inf\{W_j^{\overline{\mu}} (x): x
\in F_j\}$. Hence, for all $\mu \in \mathcal{M}_1,$
\[ \int W^{\overline{\mu}}(x) d(\mu  -
\overline{\mu})(x) = \sum_{j=1}^N \int W^{\overline{\mu}_j}(x)
d(\mu_j - \overline{\mu}_j)(x) \geq \sum_{j=1}^N
w_j^{\overline{\mu}} - w_j^{\overline{\mu}} = 0
\]
and it follows that $\overline{\mu}$ minimizes the energy
functional. With this we conclude the comments on Lemma 1.

\end{document}